\documentclass[journal,twoside,web]{ieeecolor}

\usepackage{lcsys}

\usepackage{graphics} 
\usepackage{amsmath} 
\usepackage{amssymb}  
\usepackage{amsfonts}

\usepackage{amsthm}
\usepackage{physics} 
\usepackage{mathtools} 
\usepackage{algorithm}
\usepackage{algpseudocode}
\usepackage{setspace}

\usepackage{arydshln}
\makeatletter
  \renewcommand*\env@matrix[1][*\c@MaxMatrixCols c]{%
    \hskip -\arraycolsep
    \let\@ifnextchar\new@ifnextchar
  \array{#1}}
\makeatother

\usepackage{tikz}
\usetikzlibrary{math}

\makeatletter
\let\NAT@parse\undefined
\makeatother

\usepackage{hyperref}
\hypersetup{
colorlinks=true,
linkcolor=black!80,
citecolor=darts
}

\usepackage{pgfplots}

\let\labelindent\relax
\usepackage{enumitem}

\usepackage{booktabs}
\newcommand{\behcet}{Beh\c{c}et~A\c{c}\i kme\c{s}e}

\newtheorem{thm}{Theorem}
\newtheorem{lem}{Lemma}

\newtheorem{cor}{Corollary}

\newtheorem{rem}{Remark}

\newcommand{\R}{\mathbb{R}}
\newcommand{\LL}{\mathbb{L}}
\newcommand{\E}{\mathbb{E}}
\newcommand{\K}{\mathbb{K}}
\newcommand{\D}{\mathbb{D}}

\newcommand{\pipg}{\textsc{pipg}}

\newcommand{\defeq}{\vcentcolon=}

\newcommand{\spec}{\operatorname{spec}}

\newcommand{\argmin}{\operatorname{argmin}}
\newcommand{\blkdiag}{\operatorname{blkdiag}}

\definecolor{darts}{HTML}{009670}
\definecolor{bricks}{HTML}{E74C3C}
\definecolor{fista}{HTML}{075187}
\definecolor{fistagray}{HTML}{929598}
\definecolor{steelblue}{HTML}{4682b4}
\definecolor{goldenrod}{HTML}{daa520}
\colorlet{steve}{blue!80!gray!80!green} 
\definecolor{beige}{RGB}{245,245,220}

\usetikzlibrary{shadows}
\usepackage{tcolorbox}
\tcbuselibrary{skins}

\newtcolorbox{mybox}
{
  enhanced jigsaw,
  colframe=darkgray!75,
  colback=beige!25,
  boxsep=0pt,
  drop shadow=darkgray!75!white,
  boxrule=0.75pt,
  hbox
}

\newcommand{\pipgtitle}{{\normalfont{\footnotesize{PIPG}}}}

\pgfplotsset{compat=1.17}
\usepgfplotslibrary{groupplots}
\title{Optimal Preconditioning for Online Quadratic Cone Programming}

\author{Abhinav G.\ Kamath, Purnanand Elango, and \behcet{}
\thanks{Abhinav G.\ Kamath and \behcet{} are with the William E.\ Boeing Department of Aeronautics \& Astronautics, University of Washington, Seattle, WA 98195, USA; \textbf{\texttt{\{agkamath, behcet\}@uw.edu}}.}%
\thanks{Purnanand Elango is with Mitsubishi Electric Research Laboratories, Cambridge, MA 02139, USA; \textbf{\texttt{elango@merl.com}}.}%
}

\begin{document}

\clearpage\maketitle
\thispagestyle{empty}

\begin{abstract}
First-order conic optimization solvers are sensitive to problem conditioning and typically perform poorly in the face of ill-conditioned problem data. To mitigate this, we propose an approach to preconditioning—the hypersphere preconditioner—for a class of quadratic cone programs (QCPs), i.e., conic optimization problems with a quadratic objective function, wherein the objective function is strongly convex and possesses a certain structure. This approach lends itself to factorization-free, customizable, first-order conic optimization for online applications wherein the solver is called repeatedly to solve problems of the same size/structure, but with changing problem data. We demonstrate the efficacy of our approach on numerical convex and nonconvex trajectory optimization examples, using a first-order conic optimizer under the hood.
\end{abstract}


\section{Introduction}

\IEEEPARstart{W}{e} consider the following QCP \cite{vandenberghe2010cvxopt} template:
\begin{subequations}
\begin{align}
    \underset{\xi}{\mathrm{minimize}}\quad &\frac{1}{2}\,\xi^{\top} P\,\xi + p^{\top} \xi \\
    \mathrm{subject~to}\quad &G\,\xi - g \in \LL \label{eq:cone_L}\\
    &\xi \in \E
\end{align}
\label{prob:template}%
\end{subequations}
where $\LL$ is a closed convex cone, $\E$ is a closed convex set, and $P$ is a positive definite matrix ($P \succ 0$), i.e., the objective function is strongly convex. The cone $\LL$ is a Cartesian product of closed convex cones, such as the zero cone, the nonnegative orthant cone, second-order cones (SOCs), and the cone of positive semidefinite (PSD) matrices. The set $\E$ is a Cartesian product of separable closed convex sets—such as halfspaces, boxes, $\ell_{2}$-norm balls, SOCs, etc—that possess closed-form (or easy-to-evaluate) projection operations \cite{bauschke2017convex, bauschke2018projecting}.

 The matrix $P$ is assumed to possess a structure that complies with the following requirements:
 \begin{subequations}
 \begin{gather}
 z \in \D \iff \xi \in \E \label{eq:req1} \\
 (H\,z - h \in \K) \iff (G\,\xi - g \in \LL) \label{eq:req2}
 \end{gather}
 \label{eq:reqs}%
 \end{subequations}
where $R$ is the upper-triangular Cholesky factor of $P$, i.e., $R^{\top} R = P$, $z \defeq R\,\xi$, $\D \defeq R\,\E$, $\begin{bmatrix}H & h\end{bmatrix} \defeq E\,\begin{bmatrix}G\,R^{-1} & g\end{bmatrix}$, and $\K \defeq E\,\LL$, for some positive definite diagonal matrix, $E$. Here, $S\,\mathbb{F} \defeq \{z \mid S^{-1}z \in \mathbb{F}\}$, where $S \in \R^{n \times n}$ and $\mathbb{F}$ is a convex set or cone, i.e., $z 
\in S\,\mathbb{F} \iff S^{-1}z \in \mathbb{F}$, and $\begin{bmatrix}A & b\end{bmatrix}$ denotes concatentation, where $A \in \R^{m \times n}$ and $b \in \R^{m}$.

Further, to enable online \emph{factorization-free} implementations with the elements of $P$ changing between successive calls to the solver, $P$ would need to be diagonal or have a specialized block-diagonal structure—such as the one described in \cite{kamath2023customized}—such that its Cholesky factorization is representable in closed-form. In cases wherein the Hessian of the objective function is static and only the constraint matrix is dynamic, the Cholesky factorization of $P$ (regardless of whether it is representable in closed-form) can be performed offline, such that the online component is factorization-free.

In practice, the template given by Problem \ref{prob:template} accounts for the general class of QCPs—including, but not limited to, quadratic programs (QPs), second-order cone programs (SOCPs), and semidefinite programs (SDPs)—with strongly convex quadratic objective functions, subject to the restrictions imposed on $P$.

This class of QCPs appears often in \emph{online} applications wherein optimization problems need to be solved in real-time, in a sequential fashion, with the problem size and structure remaining the same, but with the problem data changing between successive calls to the solver. Such a mode of operation (with the objective function being strongly convex) is observed most notably in model predictive control (MPC) \cite{yu2021proportional} and nonconvex trajectory optimization using sequential convex programming (SCP)—specifically a variant of SCP called sequential conic optimization (SeCO) \cite{kamath2023customized, kamath2023real}. Further, there exist online applications that involve the solution to a single convex problem (with a strongly convex objective function) \cite{elango2022customized}, which would fit well into such a framework if used in conjunction with a line search \cite{acikmese2007convex}. We note that several online applications, wherein the objective function is not necessarily strongly convex, exist as well \cite{wang2010fast, patel2011trajectory, tassa2012synthesis, liu2017real, szmuk2020successive, reynolds2020dual, yu2023real, elango2024successive}, but in this work, we restrict our focus to applications with strongly convex objective functions. While the theoretical results in this paper only hold for problems with strongly convex objective functions, problems with ill-conditioned non-strongly convex objective functions might also benefit from the proposed preconditioning procedure in practice.

First-order conic optimization solvers are attractive for: (i) real-time applications, since they are (can be) fast, (ii) implementation onboard resource-constrained embedded systems, since they only require a small code footprint, with low memory requirements in terms of both allocations and the compiled binary executable, (iii) easy \emph{customization}, which refers to exploitation of the sparsity structure of the optimization problem being solved—in a manner that enables low-dimensional matrix-vector multiplications and other dense linear algebra operations with devectorized variables—thus eliminating the need for sparse linear algebra operations entirely \cite{kamath2023customized}, (iv) easy verification and validation (V\&V), owing to the fact that they only rely on simple linear algebra operations, (v) usage within sequential convex programming algorithms for nonconvex optimization, given their amenability to warm-starting, and (vi) large-scale problems, since their performance scales well with problem size \cite{chambolle2010first, beck2019first, odonoghue2016conic, stellato2020osqp, yu2021proportional, yu2022proportional, yu2023proportional}. The performance of first-order methods, however, is highly dependent on conditioning of the problem data \cite{giselsson2014preconditioning, stellato2020osqp}—more so than solvers based on interior point methods (IPMs) that utilize second-order (Hessian) information.

Preconditioning is an operation that transforms a given matrix into another matrix of the same size but with a smaller condition number—it is a heuristic often employed to improve the performance of iterative algorithms. Simply put, applying the algorithm in question to a transformed matrix with a smaller condition number typically leads to better performance in practice (fewer iterations to convergence, for example) \cite{wathen2015preconditioning}. Several preconditioning techniques exist in the literature for various iterative numerical methods \cite{benzi2002preconditioning}.

In the context of first-order methods, diagonal preconditioners are popularly used \cite{pock2011diagonal, giselsson2014diagonal}. The computation of an exact diagonal preconditioner can be posed as a (computationally expensive-to-solve) semidefinite program \cite{giselsson2015metric, boyd1994linear}. However, in practice, matrix equilibration \cite{giselsson2015metric, sinkhorn1967concerning}—which refers to transforming the matrix in question such that the columns have equal norms and the rows have equal norms—and specifically, Ruiz equilibration \cite{ruiz2001scaling}, is a popularly used \emph{iterative} algorithm for approximate preconditioning in first-order methods that is observed to work well, while also being computationally efficient \cite{odonoghue2016conic, stellato2020osqp}. Recently, a QR factorization-based constraint matrix preconditioner was also proposed \cite{chari2024constraint}, which was shown to accelerate convergence of first-order methods. However, it requires \emph{matrix factorizations} and turns sparse problems dense, thus increasing the per-iteration cost of the solver and limiting customizability. Further, it only applies to QCPs in which $\LL$ (in Equation \ref{eq:cone_L}) only contains the zero cone. Note that all the aforementioned approaches to preconditioning are either iterative or based on (explicit) matrix factorizations.

\subsection{Contributions}

In this work, we propose a three-step procedure for preconditioning conic optimization problems that is amenable to a (mostly) analytical, factorization-free, and customizable implementation:
\begin{enumerate}[left=0pt]
    \item Minimizing the condition number of the Hessian of the objective function using hypersphere preconditioning: this step is optimal, in that the condition number of the resulting objective function matrix is unity. This was first introduced in \cite{kamath2023customized}, but it was only applied to problems where $\LL$ only contains the zero cone.
    \item Performing a (block) row-normalization of the constraint matrix, which is a simple, customization-friendly heuristic to make the constraint matrix better conditioned: although not guaranteed to reduce the condition number of the constraint matrix, this procedure has been observed to work well in practice. When used in tandem with a primal-dual conic optimization algorithm, it also provides the added benefit of scaling the dual variables favorably.
    \item Minimizing the condition number of the resulting Karush-Kuhn-Tucker (KKT) matrix by optimally scaling the objective function: in \cite{kamath2023customized}, the objective function scaling factor was manually tuned to obtain good performance, whereas in this work, we obtain a closed-form solution for the objective function scaling factor that minimizes the condition number of the KKT matrix—this step is optimal in that sense. A customizable numerical algorithm, the shifted power iteration method \cite{wilkinson1988algebraic}, can be used to efficiently estimate this value—this is the only iterative component in the entire preconditioning process. Note that this step is important, as the condition number of the KKT matrix is sensitive to the scaling of the objective function—a ``bad'' scaling can lead to poor convergence behavior, even if the first two steps are followed.
\end{enumerate}
We refer to this three-step procedure as \emph{hypersphere preconditioning}. While the end goal of preconditioning is to minimize the condition number of the KKT matrix of Problem \ref{prob:template} itself, we propose the aforementioned three-step procedure to approximately achieve that goal while ensuring that the structure of the preconditioned problem is preserved, thus enabling customization and facilitating online optimization.

\subsection{Organization}

This paper is organized as follows: Section \ref{sec:preconditioning} describes the three-step preconditioning procedure that we propose, Section \ref{sec:numerical} provides a comparison of our method against two state-of-the-art preconditioners, modified Ruiz equilibration \cite{stellato2020osqp} and the QR preconditioner \cite{chari2024constraint}, and Section \ref{sec:conclusion} concludes the paper and provides some potential avenues for future work.
\section{Preconditioning} \label{sec:preconditioning}
\subsection{Objective function: hypersphere preconditioning}
The objective function can be preconditioned using the hypersphere preconditioner \cite{kamath2023customized}, which uses the Cholesky factorization of $P$, i.e., $R^{\top}R = P$, and scales the objective function with a scalar parameter $\lambda > 0$, to transform Problem \ref{prob:template} into the following problem:
\begin{subequations}
\begin{align}
    \underset{z}{\mathrm{minimize}}\quad &\frac{\lambda}{2}z^{\top}z + \lambda\,q^{\top}z \\
    \mathrm{subject~to}\quad &\hat{G}\,z - g \in \LL \label{eq:pre_precond} \\
    &z \in \D
\end{align}
\label{prob:hypersphere}%
\end{subequations}
where $z \defeq R\,\xi$, $q \defeq R^{-\top}p$, $\hat{G} \defeq G\,R^{-1}$, and $\D \defeq R\,\E$.
This preconditioner is optimal in the sense of minimizing the condition number of the Hessian of the objective function, i.e., the condition number of the resulting objective function matrix, $\lambda\,I$, is unity. Further, its level sets are hyperspheres; hence the name. For ease of reference, we overload the term \emph{hypersphere preconditioner} to refer to the entire three-step preconditioning procedure described herewithin.
\subsection{Constraints: block row-normalization}
The constraint matrix, $\hat{G}$, can be preconditioned using row normalization, i.e., dividing each row of Equation \eqref{eq:pre_precond} by the norm of the corresponding row of $\hat{G}$, such that the rows of the preconditioned matrix have unit norms. In practice, this procedure often helps reduce the condition number of the constraint matrix, with the added benefit of not requiring transformation of the primal variable, $z$ (since the constraint matrix is only left-multiplied).

Further, row normalization should be carried out in the \emph{block} sense, i.e., by normalizing the rows of $\hat{G}$ such that the requirement given by Equation \eqref{eq:req2} holds. See \cite[Section 5]{odonoghue2016conic} for a description of this requirement. If $\LL$ only contains linear (in)equalities, such as the zero cone or the nonnegative orthant cone, then the diagonal matrix, $E \succ 0$, in $\begin{bmatrix}H & h\end{bmatrix} = E\,\begin{bmatrix} \hat{G} & g \end{bmatrix}$, can be left unrestricted. Otherwise, for each separable\footnote{The separable convex sets in $\LL$ form a partition, i.e., each set corresponds to unique components of the stacked decision variable vector, $z$.} convex set in $\LL$, only the maximum magnitude element—among the rows of $\hat{G}$ that correspond to a separable set in $\LL$—is considered in the corresponding rows of $E$.

To illustrate this, consider the following SOC constraint, where $A \in \R^{2 \times 2} \succ 0$ (diagonal), $x \in \R^{2}$, and $t \in \R$:
\begin{align}
    \norm{A\,x}_{2} \le t,~\text{i.e.},~\tilde{A}\,\tilde{x} \in \LL^{2}_{\textsc{soc}}
\end{align}
where $\norm{\cdot}_{2}$ is the Euclidean norm, $\tilde{A} \defeq \blkdiag\{A,\,1\}$, $\tilde{x} \defeq (x,\, t)$, and $\LL^{2}_{\textsc{soc}} \defeq \{(z_{1},\,z_{2},\,z_{3}) \in \R^{3} \mid \norm{(z_{1},\,z_{2})}_{2} \leq z_{3}\}$ is the second-order cone. Here, $(a,\,b) \in \R^{m + n}$, where $a \in \R^{n}$ and $b \in \R^{m}$, denotes vector concatenation. The requirement given by Equation \eqref{eq:req2} implies that a suitable choice for $E_{\textsc{soc}}$, given $E_{\textsc{soc}}\,\tilde{A}\,\tilde{x} \in E_{\textsc{soc}}\,\LL^{2}_{\textsc{soc}}$, is $c\,I_{3}$, where $c \in \R_{++}$ and $I_{3}$ is the identity matrix in $\R^{3 \times 3}$. In block row-normalization, $c \defeq \frac{1}{\max\{\abs{A_{11}},\,\abs{A_{22}},\,1\}}$, where $A_{ij}$ is the element of $A$ in the $i^{\text{th}}$ row and the $j^{\text{th}}$ column.
\subsection{Optimal objective function scaling factor}
The preconditioned problem is given as follows:
\begin{subequations}
\begin{align}
    \underset{z}{\mathrm{minimize}}\quad &\frac{\lambda}{2}z^{\top}\,z + \lambda\,q^{\top}z \\
    \mathrm{subject~to}\quad &H\,z - h \in \K \label{eq:post_precond} \\
    &z \in \D
\end{align}
\label{prob:preconditioned}%
\end{subequations}
where $\lambda > 0$ is the objective function scaling factor, and Equation \eqref{eq:post_precond} is obtained by block row-normalizing Equation \eqref{eq:pre_precond}, i.e., $\begin{bmatrix}H & h\end{bmatrix} = E\begin{bmatrix}\hat{G} & g\end{bmatrix}$, $E \succ 0$ being the diagonal matrix that performs block row-normalization and $\K \defeq E\,\LL$ being the corresponding cone.

The following lemma and corollary relate to the KKT matrix of Problem \ref{prob:preconditioned}, given by $K \defeq \begin{bmatrix}\lambda\,I & H^{\top} \\ H & 0 \end{bmatrix}$.
\begin{lem}
    For a given $\lambda > 0$ and $H \in \R^{m \times n}$, where $n > m$ and $\rank H = m$, the spectrum of $K$ is:
    \begin{align}
        \spec K =
        \left\{\lambda,\ \tfrac{\lambda \pm \sqrt{\lambda^{2} + 4\,\sigma_{1}}}{2},\ \ldots,\ \tfrac{\lambda \pm \sqrt{\lambda^{2} + 4\,\sigma_{m}}}{2}\right\} \label{eq:spectrum}
    \end{align}
    where $\sigma_{1},\,\ldots,\,\sigma_{m}$ are the squares of the singular values of $H$. The eigenvalue $\lambda$ is repeated $n-m$ times, and there are a total of $n + m$ eigenvalues.
\end{lem}
\begin{proof}
    Since $\rank H = m$ and $\lambda\,I \succ 0$, the matrix $K$ is nonsingular \cite[Section 10.1.1]{Boyd2004}.
    Consider a nontrivial eigenvector of $K$, $p \defeq (u, v) \neq 0$, where $u \in \R^{n}$ and $v \in \R^{m}$. Let $\theta$ be the corresponding eigenvalue. Then, we have $K\,p = \theta\,p$, which leads to the following equations:
    \begin{align}
        & \lambda\,u + H^{\top}v = \theta\,u \implies H^{\top}v = (\theta - \lambda)\,u \label{eq:block_row_1} \\
        & H\,u = \theta\,v \label{eq:block_row_2}
    \end{align}
    Note that $H^{\top}$ defines a one-to-one transformation, since it is full column rank. Suppose $v = 0$. This implies that $u \in \ker H$ and $(\theta - \lambda)\,u = 0$, which in turn either implies that $\theta = \lambda$ or that $u = 0$. However, since $u$ cannot be trivial, we have $\theta = \lambda$, i.e., $\theta = \lambda$ is an eigenvalue of $K$, the corresponding eigenvector being $(u, 0)$, where $u \in \ker H$. Further, there are $n-m$ eigenvectors for $\theta = \lambda$, since $\dim \ker H = n-m$.
    
    Now, multiplying Equation \eqref{eq:block_row_1} on the left by $H$ and substituting Equation \eqref{eq:block_row_2}, we get:
    \begin{align}
        H\,H^{\top}v &= (\theta - \lambda)\,H\,u = \theta\,(\theta - \lambda)\,v \label{eq:AAT}
    \end{align}
    Since $H\,H^{\top}$ is symmetric positive definite, its singular values are equal to its eigenvalues (all positive). Consider the eigenvalues of $H\,H^{\top}$, $\sigma_{k} > 0$, $k = 1,\,\ldots,\,m$. Note that these $m$ eigenvalues of $H\,H^{\top}$ have $m$ distinct corresponding eigenvectors. Now, from Equation \eqref{eq:AAT}, for each $\sigma_{k}$, we have:
    \begin{align}
        \theta\,(\theta - \lambda) = \sigma_{k} \implies \theta^{2} - \lambda\,\theta - \sigma_{k} = 0
    \end{align}
    Solving for $\theta$, we get $2\,m$ eigenvalues of $K$:
    \begin{align}
    \theta_{\pm}(\lambda, \sigma_{k}) &\defeq \tfrac{\lambda \pm \sqrt{\lambda^{2} + 4\,\sigma_{k}}}{2},\enskip\ k = 1,\,\ldots,\,m
    \end{align}

    Therefore, the spectrum of $K$ is given by Equation \eqref{eq:spectrum}, where the algebraic multiplicity of $\lambda$ is $n - m$, and the total number of eigenvalues of $K$ is $n + m$. \qedhere
\end{proof}

\begin{cor}
    For a given $\lambda > 0$ and $H \in \R^{m \times n}$, where $n > m$ and $\rank\,H = m$, the condition number of $K$ can be given by:
    \begin{align}
        \kappa(\lambda) \defeq \frac{\tfrac{\lambda + \sqrt{\lambda^2 + 4\,\sigma_{\max}}}{2}}{\min\!\left\{\lambda, \tfrac{\sqrt{\lambda^2 + 4\,\sigma_{\min}} - \lambda}{2}\right\}}
        \label{eq:kappa}
    \end{align}
    where $\sigma_{\max}$ and $\sigma_{\min}$ are the squares of the largest and smallest singular values of $H$, respectively.
\end{cor}
\begin{proof}
    Since $\theta_{+}(\lambda, \sigma_{i}) \geq \max\{\lambda,\,|\theta_{-}(\lambda, \sigma_{j})|\}$, $1 \leq i,\,j \leq m$, the condition number of $K$ as a function of $\lambda$ can be given by $\kappa(\lambda) = \frac{\delta_{\max}(\lambda)}{\delta_{\min}(\lambda)}$, where:
    {\scriptsize
    \begin{flalign}
        \hspace*{-1em} \delta_{\max}(\lambda) &\defeq \max_{k} \theta_{+}(\lambda, \sigma_{k}) = \max_{k}\tfrac{\lambda + \sqrt{\lambda^2 + 4\,\sigma_{k}}}{2}\label{eq:max_singular_value} \\ 
        \hspace*{-1em} \delta_{\min}(\lambda) &\defeq \min\{\lambda,\,\min_{k}|\theta_{-}(\lambda, \sigma_{k})|\} = \min\!\left\{\lambda,\,\min_{k} \tfrac{\sqrt{\lambda^2 + 4\,\sigma_{k}} - \lambda}{2}\right\} \label{eq:min_singular_value}
    \end{flalign}
    }%
    The term $\theta_{+}(\lambda, \sigma_{k})$ in Equation \eqref{eq:max_singular_value} attains its maximum when $\sigma_{k} = \sigma_{\max}$.
    \begin{align}
        \therefore\ \delta_{\max}(\lambda) = \tfrac{\lambda + \sqrt{\lambda^2 + 4\,\sigma_{\max}}}{2}
    \end{align}
    The term $|\theta_{-}(\lambda, \sigma_{k})|$ in Equation \eqref{eq:min_singular_value} attains its minimum when $\sigma_{k} = \sigma_{\min}$.
    \begin{align}
        \therefore\ \delta_{\min}(\lambda) = \min\!\left\{\lambda, \tfrac{\sqrt{\lambda^2 + 4\,\sigma_{\min}} - \lambda}{2}\right\} \label{eq:min_singular_value_attained}
    \end{align}
    \vspace{-1em}
    \begin{align*}
        \therefore\ \kappa(\lambda) = \tfrac{\tfrac{\lambda + \sqrt{\lambda^2 + 4\,\sigma_{\max}}}{2}}{\min\left\{\lambda, \tfrac{\sqrt{\lambda^2 + 4\,\sigma_{\min}} - \lambda}{2}\right\}} \tag*{\qedhere}
    \end{align*}
\end{proof}
\begin{thm}
    For a given $H \in \R^{m \times n}$, where $n > m$ and $\rank\,H = m$,
    \vspace{-1em}
    \begin{align}
    \lambda^{\star} \defeq \underset{\lambda}{\argmin}\,\kappa(\lambda) = \sqrt{\frac{\sigma_{\min}}{2}} \label{eq:lambda_star}
    \end{align}
    where $\lambda > 0$, $\sigma_{\min}$ is the square of the smallest singular value of $H$, and $\kappa(\lambda)$ is given by Equation \eqref{eq:kappa}. \label{thm:min_cond}
\end{thm}
\begin{proof}
    Let $f_{1}(\lambda) \defeq \lambda$ and $f_{2}(\lambda) \defeq \tfrac{\sqrt{\lambda^{2} + 4\,\sigma_{\min}} - \lambda}{2}$.
    From Equation \eqref{eq:min_singular_value_attained}, $\delta_{\min}(\lambda) = \min\!\left\{f_{1}(\lambda),\,f_{2}(\lambda)\right\}$. It is clear that $f_{1}$ is a positive and strictly increasing function of $\lambda$. Now, taking the derivative of $f_{2}(\lambda)$ with respect to $\lambda$, we get:
    \begin{align}
        \tfrac{df_{2}(\lambda)}{d\lambda} = \tfrac{1}{2}\!\left(\tfrac{\lambda}{\sqrt{\lambda^2 + 4\,\sigma_{\min}}} - 1\right)
    \end{align}
    Since $\sqrt{\lambda^{2} + 4\,\sigma_{\min}} > \lambda$, we conclude that $\tfrac{df_{2}(\lambda)}{d\lambda} < 0$, and hence, $f_{2}$ is a positive and strictly decreasing function of $\lambda$.
    \begin{align}
        \therefore\ \min\!\left\{f_{1}(\lambda),\,f_{2}(\lambda)\right\} = \tfrac{1}{\max\left\{\tfrac{1}{f_{1}(\lambda)},\,\tfrac{1}{f_{2}(\lambda)}\right\}}
    \end{align}
    Now, let $f_{0}(\lambda) \defeq \tfrac{\lambda + \sqrt{\lambda^2 + 4\,\sigma_{\max}}}{2}$. We see that $\tfrac{df_{0}(\lambda)}{d\lambda} > 0$, and hence, $f_{0}$ is a positive and strictly increasing function of $\lambda$.
    {\footnotesize
    \begin{equation}
        \therefore\ \kappa(\lambda) = f_{0}(\lambda)\,\max\!\left\{\tfrac{1}{f_{1}(\lambda)},\,\tfrac{1}{f_{2}(\lambda)}\right\} = \max\!\left\{\tfrac{f_{0}(\lambda)}{f_{1}(\lambda)},\,\tfrac{f_{0}(\lambda)}{f_{2}(\lambda)}\right\} \label{eq:kappa_max_form}
    \end{equation}
    \begin{equation}
        \therefore\ \min_{\lambda}\,\kappa(\lambda) = \min_{\lambda}\,\max\!\left\{\tfrac{f_{0}(\lambda)}{f_{1}(\lambda)},\,\tfrac{f_{0}(\lambda)}{f_{2}(\lambda)}\right\}
    \end{equation}
    }%
    Since $\tfrac{f_{0}(\lambda)}{f_{1}(\lambda)}$ is strictly decreasing and $\tfrac{f_{0}(\lambda)}{f_{2}(\lambda)}$ is strictly increasing, the minimizer of Equation \eqref{eq:kappa_max_form} satisfies $\tfrac{f_{0}(\lambda)}{f_{1}(\lambda)} = \tfrac{f_{0}(\lambda)}{f_{2}(\lambda)}$, i.e., $f_{1}(\lambda) = f_{2}(\lambda)$ (since $f_{0}(\lambda)$ is positive), as depicted in Fig.\ \ref{fig:min_max}. $\enskip\ \therefore\,\lambda^{\star} = \tfrac{\sqrt{\lambda^{\star^{2}} + 4\,\sigma_{\min}} - \lambda^{\star}}{2} \implies \lambda^{\star} = \sqrt{\tfrac{\sigma_{\min}}{2}}$.
\end{proof}
\vspace{-1.25em}
\begin{figure}[H]
    \centering
    {
    \resizebox{0.875\linewidth}{!}{\begin{tikzpicture}

\footnotesize

    \def\fA{0.0125 * (\x + 1.125)^3 + 1}  
    \def\fB{3 - 0.5 * \x}  
    
    \draw[line width = 1,line cap=round,->] (-0.5,0) -- (5.5,0) node[right] {$x$};
    \draw[line width = 1,line cap=round,->] (0,-0.5) -- (0,4);
  
    \draw[line width = 2,line cap=round,domain=-0.5:4.5,smooth,variable=\x,darts!75] plot ({\x},{\fA}) node[right] {$f(x)$};
    \draw[line width = 2,line cap=round,domain=-0.5:4.5,smooth,variable=\x,bricks!75] plot ({\x},{\fB}) node[right] {$g(x)$};

    \tikzmath{\xxx = 2.65246; \yyy = 1.67377;} 
  
    \coordinate (I) at (\xxx, \yyy);
  
    \draw[very thick,line cap=round,dashed,darkgray!75] (I) -- (\xxx, 0) node[below] {{\color{black}{$x^\star$}}};
    \draw[very thick,line cap=round,dashed,darkgray!75] (I) -- (0, \yyy) node[left] {{\color{black}{$f(x^\star) = g(x^\star)$}}};

    \draw[line width = 2.375,draw opacity=0.75,line cap=round,domain=-0.475:\xxx,smooth,variable=\x,white] plot ({\x},{\fB+0.065});
    \draw[line width = 2.375,draw opacity=0.75,line cap=round,domain=\xxx:4.475,smooth,variable=\x,white] plot ({\x},{\fA+0.065});
    
    \draw[line width = 1.375,line cap=round,domain=-0.475:\xxx,smooth,variable=\x,steve!87.5] plot ({\x},{\fB+0.065});
    \draw[line width = 1.375,line cap=round,domain=\xxx:4.475,smooth,variable=\x,steve!87.5] plot ({\x},{\fA+0.065});
    
    \node[left, steve!87.5] at (3.55, 2.775) {$\max\{f(x),\,g(x)\}$};
    \node[draw=white, draw opacity=0.75, minimum size=4.25pt, inner sep=0pt] at (\xxx, \yyy) {};
    \node[fill=steve!87.5, minimum size=4pt, inner sep=0pt] at (\xxx, \yyy) {};
    \fill[draw=white, draw opacity=0.75] (\xxx, \yyy) circle (1.25pt);
    \fill[darkgray!87.5] (\xxx, \yyy) circle (1.125pt);
    \draw[line width = 1,darkgray!87.5,line cap=round,->] (\xxx, \yyy) -- (\xxx + 0.75, \yyy) node[right] {{\color{black}{$\underset{x}{\min}\!\;\max\{f(x),\,g(x)\}$}}};
  
  \end{tikzpicture}}
    }
    \vspace{-0.4375em}
    \caption{The minimizer of the maximum of a strictly increasing function, $f(x)$, and a strictly decreasing function, $g(x)$, satisfies $f(x) = g(x)$.}
    \label{fig:min_max}
\end{figure}
\vspace{-1.25em}
\begin{cor}
    With the three-step preconditioning procedure described, $\kappa(\lambda^{\star}) \ge 2$, i.e., the condition number of the preconditioned KKT matrix is lower-bounded by $2$. \label{cor:cond_lower_bound}
\end{cor}
\begin{proof}
    By substituting Equation \eqref{eq:lambda_star} in Equation \eqref{eq:kappa}, we get:
    $$\kappa(\lambda^{\star}) = \tfrac{1 + \sqrt{1 + 8\,\chi}}{2}$$ 
    where $\chi \defeq \frac{\sigma_{\max}}{\sigma_{\min}}$ is the condition number of $H\,H^{\top}$. When $\chi = 1$, the lower bound on $\kappa(\lambda^{\star})$ is tight.
\end{proof}
Corollary \ref{cor:cond_lower_bound} provides a limit on how much we can minimize the condition number of the KKT matrix with the three-step preconditioning procedure described: if $H\,H^{\top}$ happens to be perfectly conditioned, the condition number of the preconditioned KKT matrix is $2$. In practice, row normalization is effective in reducing the condition number of $H\,H^{\top}$, often by a few orders of magnitude, although rarely to unity. So, although a condition number of unity for the KKT matrix is unattainable (due to being lower-bounded by $2$), and despite its lower bound usually not being tight (as a result of $H\,H^{\top}$ not being perfectly conditioned), the proposed preconditioning procedure can reduce the condition number of the KKT matrix enough to significantly improve the performance of first-order QCP solvers, without sacrificing any of the following features: being (i) mostly analytical (the only iterative component is shifted power iteration), (ii) entirely factorization-free, and (iii) amenable to customization—all of which are beneficial for online applications.
\vspace{-0.5em}
\section{Numerical results} \label{sec:numerical}
We first consider the convex optimal control problem from \cite{yu2021proportional} to demonstrate the efficacy of the hypersphere preconditioner. To enable comparison of our method with modified Ruiz equilibration \cite{stellato2020osqp}, we replace the $2$-norm ball constraints with their $\infty$-norm counterparts. Further, in compliance with the requirement given by Equation \eqref{eq:req1}, we choose the following state cost matrix: $Q_{t} \defeq Q = \operatorname{diag}\{1, 1, 0.5, 0.5\},\ t = 1, \ldots, T-1$, where $T$ is chosen to be $50$. We scale the terminal state cost matrix by a scaling factor, $\gamma > 0$, i.e., $Q_{T} \defeq \gamma\,Q$. By progressively increasing $\gamma$, we obtain smaller and smaller terminal state tracking errors, but at the cost of increasing ill-conditioning in the problem. 

With that, we compare the hypersphere preconditioner with two state-of-the-art preconditioning techniques: modified Ruiz equilibration \cite{stellato2020osqp} and the QR preconditioner \cite{chari2024constraint}, in terms of (i) minimizing the condition number of the KKT matrix of the preconditioned problem, and (ii) performance of the proportional-integral projected gradient method (\pipg{}) \cite{yu2021proportional, yu2021infeasibility, yu2022proportional, yu2023proportional}, a first-order primal-dual method for conic optimization, that has recently gained popularity for problems that fit the trajectory optimization template \cite{yu2023real, loureno2023verification, kamath2023real, kamath2023customized, chari2024fast, chari2024spacecraft, mangalore2024neuromorphic, doll2025hardware}.

\begin{figure}[H]
    \centering
    \includegraphics[width=0.95\linewidth]{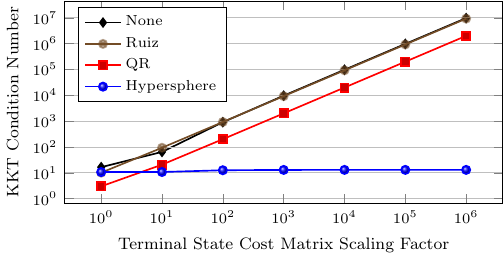}
    \vspace{-0.5em}
    \caption{The condition number of the preconditioned KKT matrix as a function of the terminal state cost matrix scaling factor, $\gamma$. ``None'' refers to the case where no preconditioning is applied.}
    \label{fig:cond}
\end{figure}

\begin{figure}[H]
    \centering
    \includegraphics[width=0.95\linewidth]{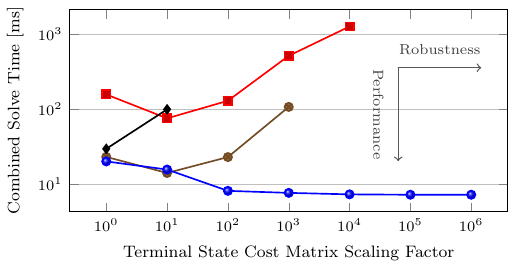}
    \vspace{-0.5em}
    \caption{The combined solve time (preconditioning and convex solve) as a function of $\gamma$. The combined solve time is only plotted if \pipg{} was able to converge (within $10^{5}$ iterations).}
    \label{fig:solve-time}
\end{figure}

All the preconditioners and the solver are implemented in C, the code for which is generated using the \textsc{matlab} Coder. All the problems are first solved to high-accuracy using an off-the-shelf interior-point method (constituting the ``ground truth''), and \pipg{} is deemed to have converged if and when the relative error between the \pipg{} solution and the ground truth drops below $0.5\%$. The code is available at {\footnotesize\url{https://github.com/UW-ACL/optimal-preconditioning}}. Fig.\ \ref{fig:cond} shows the effect of the preconditioners on the condition number of the preconditioned KKT matrix, and Fig.\ \ref{fig:solve-time} and Table \ref{tab:PIPG-iters} show their effect on the overall solver performance.

 To estimate $\sigma_{\min}$ within the hypersphere preconditioner, we use the shifted power iteration method \cite{wilkinson1988algebraic}, the performance statistics for which are provided in Table \ref{tab:shifted-power}. The presolve time includes the preconditioning time and the time taken by the regular power iteration method to estimate $\sigma_{\max}$, which \pipg{} uses for step-size determination; however, the presolve step is largely dominated by shifted power iteration. The absolute and relative tolerances for the termination of shifted power iteration are each set to $10^{-9}$. The convergence rate of shifted power iteration is proportional to the ratio $\frac{\sigma_{\max} - \sigma_{\min_{2}}}{\sigma_{\max_{\phantom{1}}}\!\!\!- \sigma_{\min}}$, where $\sigma_{\min_{2}}$ is the second smallest eigenvalue of $H\,H^{\top}$. We observe that shifted power iteration suffers from slow convergence when the smallest eigenvalues of $H\,H^{\top}$ are clustered (explaining the slower convergence for smaller values of $\gamma$)—this is a limitation of the method.
\vspace{-1.5em}
\vspace{-1em}
\begin{table}[H]
\scriptsize
    \centering
    \begin{tabular}{cccccccc}
\toprule
\textbf{None}        & $3081$  & $11209$ &  ${\color{bricks}{10^{5}}}$  &  ${\color{bricks}{10^{5}}}$   &  ${\color{bricks}{10^{5}}}$   &  ${\color{bricks}{10^{5}}}$   &  ${\color{bricks}{10^{5}}}$   \\
\textbf{Ruiz}        & $1101$  & ${\color{darts}{\mathbf{237}}}$   & $1262$  & $10774$ &  ${\color{bricks}{10^{5}}}$   &    ${\color{bricks}{10^{5}}}$ &  ${\color{bricks}{10^{5}}}$   \\
\textbf{QR}          & $2240$  & $2098$  & $2014$  & $3894$  & $28588$ &  ${\color{bricks}{10^{5}}}$   & ${\color{bricks}{10^{5}}}$    \\
\textbf{Hypersphere} & ${\color{darts}{\mathbf{684}}}$   & $672$   & ${\color{darts}{\mathbf{558}}}$   & ${\color{darts}{\mathbf{529}}}$   & ${\color{darts}{\mathbf{525}}}$   & ${\color{darts}{\mathbf{524}}}$   & ${\color{darts}{\mathbf{524}}}$  \\
\midrule
$\gamma$ & $1$ & $10$ & $10^{2}$ & $10^{3}$ & $10^{4}$ & $10^{5}$ & $10^{6}$ \\
\bottomrule
\end{tabular}
    \caption{The number of \pipg{} iterations as a function of $\gamma$. The {\color{darts}{green}} values indicate the least number of \pipg{} iterations to convergence. The {\color{bricks}{red}} values indicate that \pipg{} hit the maximum number of iterations ($10^{5}$) and failed to converge.}
    \label{tab:PIPG-iters}
\end{table}
\vspace{-1.75em}
\begin{table}[H]
\scriptsize
    \centering
    \begin{tabular}{cccccccc}
\toprule
\textbf{Iterations}        & $6081$  & $4758$ &  $1413$  &  $1304$   &  $1294$   &  $1293$   &  $1293$   \\
$\frac{\textbf{PT}}{\textbf{CST}} (\%)$   & $68.06$  & $62.74$   & $39.2$  & $37.37$ &  $36.91$   &   $36.46$ &  $36.64$   \\
\midrule
$\gamma$ & $1$ & $10$ & $10^{2}$ & $10^{3}$ & $10^{4}$ & $10^{5}$ & $10^{6}$ \\
\bottomrule
\end{tabular}
    \caption{Performance of shifted power iteration within hypersphere preconditioning. \textbf{pt}: presolve time; \textbf{cst}: combined solve time.}
    \label{tab:shifted-power}
\end{table}
\vspace{-1em}
Hypersphere preconditioning was designed with the express goal of mitigating ill-conditioning in the objective function Hessian \cite{kamath2023customized}, and we observe that it performs the best, with increasing ill-conditioning, across all metrics. Although \pipg{} requires longer solve times with the QR preconditioner owing to the loss of sparsity, the QR preconditioner could prove beneficial in solving problems that have dense and ill-conditioned constraint matrices. Finally, modified Ruiz equilibration is a balanced preconditioning approach, and is the only one among the three that does not mandate strong convexity and explicitly accounts for the affine term in the objective function.

Next, to demonstrate the benefit of choosing the optimal objective function scaling factor, we consider a numerical example involving a practical online application. More specifically, we consider the nonconvex multi-phase rocket landing guidance problem from \cite{kamath2023real}, which is solved using sequential conic optimization (SeCO). This involves solving a sequence of QCP subproblems that have strongly convex objective functions, using \pipg{}, with the application of the hypersphere preconditioner.

\begin{figure}[H]
    \centering
    \includegraphics[width=0.95\linewidth]{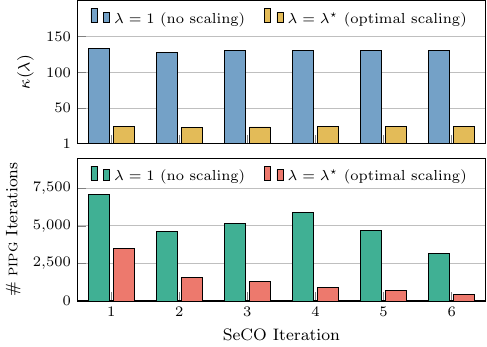}
    \vspace{-0.5em}
    \caption{The effect of the optimal objective function scaling factor, $\lambda^{\star}$—on (i) the condition number, $\kappa(\lambda)$, of the KKT matrix, and (ii) the number of \pipg{} iterations—at each iteration of SeCO for the multi-phase rocket landing guidance problem \cite{kamath2023real}.}
    \label{fig:results}
\end{figure}

The guidance problem is solved in $6$ SeCO iterations to a predetermined open-loop terminal-error accuracy for the translation states ($< 1$ m in position and $< 0.5$ m\;s$^{-1}$ in velocity). Fig.\ \ref{fig:results} shows a comparison of cases with no objective function scaling ($\lambda = 1$) and the optimal objective function scaling ($\lambda = \lambda^{\star}$). We observe a clear reduction in (i) the condition number of the KKT matrix for each subproblem, and (ii) the number of \pipg{} iterations to convergence for each of the subproblems, thus demonstrating the importance of properly scaling the cost function.
\vspace{-0.5em}
\section{Conclusion} \label{sec:conclusion}
We propose a three-step preconditioning procedure, the hypersphere preconditioner, to improve the performance of conic optimization solvers for online applications that involve solving quadratic cone programs (QCPs) with strongly convex objective functions. This preconditioner is amenable to a mostly analytical, fully factorization-free, and customizable implementation, while significantly improving the performance of first-order conic optimization solvers, making it particularly beneficial for online applications that involve solving a sequence of QCPs with dynamically changing problem data. We derive an analytical expression for the optimal objective function scaling factor in the sense of minimizing the condition number of the KKT matrix of the preconditioned problem, and obtain a lower bound on it. We demonstrate the efficacy of the hypersphere preconditioner by means of numerical experiments in both convex optimization and sequential convex programming settings.

Future work involves expanding the class of problems to which the proposed preconditioning procedure can apply. Generalizing the preconditioning procedure to problems with non-strongly convex objective functions would allow for its use within most NMPC and SCP-based trajectory optimization frameworks \cite{malyuta2021advances, malyuta2022convex, mceowen2023high, chari2024fast, chari2024spacecraft, uzun2024successive, kim2024six}.
\vspace{-0.5em}
\section*{Acknowledgments}

The authors thank Govind M.\ Chari for discussions on row normalization and the QR preconditioner, and Danylo Malyuta for detailed feedback on an initial draft.
\vspace{-0.5em}

\bibliographystyle{ieeetr}
\bibliography{root}

\section*{Appendix}
\subsection{Optimal step-size ratio in {\small PIPG}}
Considering Problem \ref{prob:preconditioned}, the primal and dual step-sizes in \pipg{} \cite{yu2023proportional}, $\alpha$ and $\beta$, respectively, satisfy the conditions $\alpha,\,\beta > 0$, $\alpha (\|\lambda\,I\| + \beta\|H\|^{2}) < 1$, where $\|\square\|$ is defined to be the largest singular value of $\square$. One choice for $\alpha$ and $\beta$ is obtained by parameterizing $\beta$ in terms of $\alpha$, i.e., $\beta \defeq \omega^{2}\alpha$, $\omega \in \R_{++}$, setting the strict inequality to an equality, and solving for $\alpha$, i.e., $
    \alpha = \tfrac{2}{\lambda + \sqrt{\lambda^{2} + 4\,\omega^{2}\,\sigma_{\max}}}$ and 
    $\beta = \omega^{2}\alpha$, 
where $\sigma_{\max} \defeq \|H\|^{2}$.
The \pipg{} iterates for Problem \ref{prob:preconditioned} are given by (see \cite{yu2023proportional} for more details):
\begin{subequations}
\begin{align}
    z^{+} &= \Pi_{\D}[\zeta - \alpha\,(\lambda\,\zeta + q + H^{\top} \eta)] \\
    w^{+} &= \Pi_{\K^{\circ}}\!\!\:[\eta + \beta\,(H(2\,z^{+}-\zeta) - h)] \label{eq:dual_update} \\
    \zeta^{+} &= (1 - \rho)\,\zeta + \rho\,z^{+} \\
    \eta^{+} &= (1 - \rho)\,\eta + \rho\,w^{+} \label{eq:dual_extrapolation}
\end{align}
\end{subequations}
The step-size ratio is given by:
\begin{align}
    \frac{\beta}{\alpha} = \omega^{2}\label{eq:step_size_ratio}
\end{align}
Typically, the solver parameter, $\omega$, is manually tuned to obtain good performance \cite{yu2023proportional}, the tuning process itself being unintuitive in practice. While \cite{chari2024constraint} provides a different step-size rule than the one considered here and adopts an adaptive heuristic based on approximately minimizing the primal-dual gap, in this work, we establish the relationship between the step-size ratio and the objective function scaling factor and find the connection between their optimal values.
\begin{thm}
    Scaling the step-size ratio in \pipgtitle{} by a factor of $s^2$, $s > 0$, is equivalent to scaling the objective function of the corresponding QCP by a factor of $\frac{1}{s}$, i.e., scaling $\omega$ by $s$ is equivalent to scaling $\lambda$ by $\frac{1}{s}$.\label{thm:pipg}
\end{thm}
\begin{proof}
    Scaling $\omega$ by $s$, we get the primal-dual step-sizes $\tilde{\alpha} \defeq \tfrac{2}{\lambda + \sqrt{\lambda^{2} + 4\,s^{2}\,\omega^{2}\,\sigma_{\max}}}$ and $\tilde{\beta} \defeq s^{2}\,\omega^{2}\,\alpha$, respectively—this corresponds to scaling the step-size ratio by $s^{2}$. Now, we have $\tilde{\alpha} = \frac{\hat{\alpha}}{s}$ and $\tilde{\beta} = s\,\omega^{2}\,\hat{\alpha}$, where $\hat{\alpha} \defeq \tfrac{2\,s}{\lambda + \sqrt{\lambda^{2} + 4\,s^{2}\,\omega^{2}\,\sigma_{\max}}}$. Let $\hat{\beta} \defeq \omega^{2}\,\hat{\alpha}$. From \cite[Section II, Lemma 3, (ii)]{Ingram1991}, the projection onto a convex cone is nonnegatively homogeneous (and hence, positively homogeneous). Therefore, $\Pi_{\K^{\circ}}[c\,z] = c\,\Pi_{\K^{\circ}}[z]$, $c > 0$. Further, defining $\hat{w} \defeq \frac{w}{s}$ and $\hat{\eta} \defeq \frac{\eta}{s}$, dividing Equations \eqref{eq:dual_update} and \eqref{eq:dual_extrapolation} by $s$, and invoking \cite[Section II, Lemma 3, (ii)]{Ingram1991} on Equation \eqref{eq:dual_update}, we get:
\begin{subequations}
\begin{align}
    z^{+} &= \Pi_{\D}\!\left[\zeta - \hat{\alpha}\left(\tfrac{\lambda}{s}\,\zeta + \tfrac{\lambda}{s}\,q + H^{\top} \hat{\eta}\right)\right] \\
    \hat{w}^{+} &= \Pi_{\K^{\circ}}\!\!\:[\hat{\eta} + \hat{\beta}\,(H(2\,z^{+}-\zeta) - h)] \\
    \zeta^{+} &= (1 - \rho)\,\zeta + \rho\,z^{+} \\
    \hat{\eta}^{+} &= (1 - \rho)\,\hat{\eta} + \rho\,\hat{w}^{+}    
\end{align}
\end{subequations}%
which are the \pipg{} iterates for the following equivalent problem, Problem \ref{prob:preconditioned_reformulated}, but with step-sizes $\hat{\alpha}$ and $\hat{\beta}$:
\begin{subequations}
\begin{align}
    \underset{z}{\mathrm{minimize}}\quad &\frac{\hat{\lambda}}{2}\,z^{\top}\,z + \hat{\lambda}\,q^{\top}z \\
    \mathrm{subject~to}\quad &H\,z - h \in \K \\
    &z \in \D
\end{align}
\label{prob:preconditioned_reformulated}
\end{subequations}
where $\hat{\alpha} = \frac{2}{\hat{\lambda} + \sqrt{\hat{\lambda}^{2} + 4\,\omega^{2}\,\sigma_{\max}}}$, $\hat{\beta} = \omega^{2}\,\hat{\alpha}$, and $\hat{\lambda} \defeq \frac{\lambda}{s}$.
\end{proof}
\begin{rem}
As a consequence of Theorem \ref{thm:pipg}, tuning $\lambda$ is equivalent to tuning $\omega$, i.e., there is a one-to-one mapping between $\lambda$ and $\omega$.
\end{rem}
\begin{cor}
    For a given $\lambda > 0$ in Problem \ref{prob:preconditioned}, the optimal value for the \pipgtitle{} solver parameter, $\omega$—in terms of minimizing the condition number of $K$—is given by $\omega^{\star} \defeq \lambda\,\sqrt{\frac{2}{\sigma_{\min}}}$.\label{cor:opt_param}
\end{cor}
\begin{proof}
    This directly follows from Theorems \ref{thm:min_cond} and \ref{thm:pipg}.
\end{proof}
\subsection{Shifted power iteration}
The maximum singular value of $M \defeq H\,H^{\top} \in \R^{m \times m}$, which is a parameter that factors into the step-sizes of \pipg{}, can be efficiently estimated using the power iteration method \cite{trefethen1997numerical}, which, in turn, can be customized to the trajectory optimization template for efficient implementation \cite{kamath2023customized}. For the power iteration method to be convergent, the magnitude of the dominant eigenvalue must be strictly greater than the magnitude of every other eigenvalue. We make the assumption that the eigenvalues of matrix $M$ satisfy this condition. The optimal solver parameters in {\pipg}—given by Corollary \ref{cor:opt_param}—however, require an estimate for the \textit{minimum} singular value of $M$.

General methods to compute the minimum singular value of a matrix, such as  the inverse iteration method or the singular value decomposition (SVD), are generally computationally expensive and not amenable to customization, i.e., they are not structure-exploiting, thus making them unsuitable for real-time applications that require online computations. However, since $M \succ 0$, we can use the shifted power iteration method to estimate its smallest singular value \cite{wilkinson1988algebraic}, which is described in Algorithm \ref{alg:shifted_power}, where $\sigma_{\min}(A)$ is defined to be the smallest eigenvalue of symmetric matrix $A \succ 0$.
\begin{algorithm}
\small
\caption{Shifted power iteration to estimate $\sigma_{\min}(H\,H^{\top})$}\label{alg:shifted_power}
    \vspace{0.5em}
    \begin{flushleft}
        \textbf{Inputs:} $H$, $w$, $\epsilon_{\mathrm{abs}}$, $\epsilon_{\mathrm{rel}}$, $\epsilon_{\mathrm{buff}}$, $j_{\max}$
    \end{flushleft}
    \begin{algorithmic}[1]
    \Require $\norm{w}_{2} > 0$
    \vspace{1em}
    \State $\tilde{\sigma} \leftarrow \norm{w}_{2}$\label{line:init} \Comment{initialization}\vspace{1ex}
    \For {$j \leftarrow 1, \ldots,\, j_{\max}$}\vspace{1ex}
    \State $z \leftarrow H^{\top}w$\label{line:z}\vspace{1ex}
    \State $w \leftarrow \frac{1}{\tilde{\sigma}}\,(H\,z - \sigma_{\max}\,w)$\label{line:w} \Comment{$\sigma_{\max} \defeq \|H\|^{2}$}\vspace{1ex}
    \State $\tilde{\sigma}^{\star} \leftarrow \norm{w}_{2}$\label{line:sigma_star}\vspace{1ex}
    \If {$\abs{\tilde{\sigma}^{\star} - \tilde{\sigma}} \le \epsilon_{\mathrm{abs}} + \epsilon_{\mathrm{rel}}\,\max\{\tilde{\sigma}^{\star},\,\tilde{\sigma}\}$}
    \State \textbf{break}
    \ElsIf {$j < j_{\max}$}
    \State $\tilde{\sigma} \leftarrow \tilde{\sigma}^{\star}$
    \EndIf\vspace{1ex}
    \EndFor\vspace{1ex}
    \State $\sigma_{\min} \leftarrow (1 - \epsilon_{\mathrm{buff}})\,(\sigma_{\max} - \tilde{\sigma}^{\star})$
    \Comment{buffer the overestimate}\vspace{0.125em}
    \end{algorithmic}
    \begin{flushleft}
        \textbf{Return:} $\sigma_{\min}$ \Comment{$\approx \sigma_{\min}(H\,H^{\top})$}
    \end{flushleft}
\end{algorithm}
First, we perform a spectral shift on $M$ as follows:
\begin{align}
    \widetilde{M} \defeq M - \sigma_{\max}\,I
\end{align}
where $\sigma_{\max}$ is the largest singular value of $M$. This shift annihilates the largest eigenvalue of $M$ and consequently \textit{deflates} it to form $\widetilde{M}$. To ensure that the shifted power iteration method is convergent, we make the assumption that the magnitude of the dominant eigenvalue of $\widetilde{M}$ is strictly greater than the magnitude of every other eigenvalue. Let $v$ be the eigenvector corresponding to an arbitrary eigenvalue (singular value) of $M$, $\sigma$, i.e., $M\,v = \sigma\,v$. The matrix-vector product $\widetilde{M}\,v$ yields
\begin{align}
    \widetilde{M}\,v = M\,v - \sigma_{\max}\,I\,v = (\sigma - \sigma_{\max})\,v
\end{align}
Therefore, $\sigma - \sigma_{\max}$ is an eigenvalue of $\widetilde{M}$. Further, since $\sigma - \sigma_{\max} \le 0$, matrix $\widetilde{M}$ is symmetric negative semidefinite. Therefore, the power iteration method can be used to find the absolute value of its largest magnitude eigenvalue (largest singular value), which is nothing but $\widetilde{\sigma} \defeq |\sigma_{\min} - \sigma_{\max}| = \sigma_{\max} - \sigma_{\min}$, where $\sigma_{\min}$ is the smallest singular value of $M$. Finally, $\sigma_{\min}$ can be obtained by subtracting $\widetilde{\sigma}$ from $\sigma_{\max}$, i.e., $\sigma_{\min} = \sigma_{\max} - \widetilde{\sigma}$.

Note that Algorithm \ref{alg:shifted_power} is amenable to customization, as shown in \cite[Algorithm 5]{kamath2023customized}.

\end{document}